\documentclass[letterpaper, 10 pt, conference]{ieeeconf}  
\usepackage[utf8]{inputenc}
\usepackage{amsmath}
\usepackage{amssymb}
\usepackage{amsfonts}
\usepackage{mathtools}
\usepackage{float}
\newtheorem{lemma}{Lemma}
\newtheorem{remark}{Remark}
\newtheorem{assumption}{Assumption}
\newtheorem{theorem}{Theorem}
\newtheorem{definition}{Definition}
\newtheorem{example}{Example}
\usepackage{mathtools}  
\usepackage{xcolor}
\mathtoolsset{showonlyrefs}

\IEEEoverridecommandlockouts                              

\makeatletter
\def\endthebibliography{%
  \def\@noitemerr{\@latex@warning{Empty `thebibliography' environment}}%
  \endlist
}
\makeatother

\title{\LARGE \bf
\textcolor{black}{Regret~Minimization~in~Scalar,~Static,~Non-linear~Optimization~Problems}
}

\author{Ying Wang,  Mirko Pasquini, Kévin Colin, Håkan Hjalmarsson
\thanks{*This work was supported by VINNOVA Competence Center AdBIOPRO, contract [2016-05181]  and by the Swedish Research Council through the research environment NewLEADS (New Directions in Learning Dynamical Systems), contract [2016-06079] and by Wallenberg AI, Autonomous Systems and Software Program (WASP), funded by Knut and Alice Wallenberg Foundation.}
\thanks{The authors are with the Division of Decision and Control Systems,
	KTH Royal Institute of Technology, 10044 Stockholm, Sweden.
	Mirko Pasquini, Kévin Colin and Håkan Hjalmarsson are also with the Center for Advanced Bio Production AdBIOPRO, KTH Royal Institute of Technology, 10044 Stockholm, Sweden (e-mail:yinwang@kth.se; pasqu@kth.se; kcolin@kth.se; hjalmars@kth.se).}
}

\begin{document}
\maketitle
\thispagestyle{empty}
\pagestyle{empty}

\begin{abstract}
We study the problem of determining an \textcolor{black}{effective exploration} strategy in static and non-linear optimization problems, which depend on an unknown scalar parameter to be learned from online collected noisy data.
An optimal trade-off between exploration and exploitation is crucial for effective optimization under uncertainties, and to achieve this we consider a cumulative regret minimization approach over a finite horizon, with each time instant in the horizon characterized by a stochastic exploration signal, whose variance is to be designed. 
\textcolor{black}{We aim to extend the well-established concepts of regret minimization from linear to non-linear systems, with a focus on the subsequent conceptual differences and challenges.
Thus}, under an idealized assumption on an appropriately defined information function associated with the excitation, we are able to show that an \textcolor{black}{optimal exploration} strategy is either to use no exploration at all (called lazy exploration) or adding an exploration excitation only at the first time instant of the horizon (called immediate exploration).
A quadratic numerical example is presented to demonstrate the effectiveness of the proposed strategy.
\end{abstract}

\vspace{-0.1cm}\section{Introduction}
\label{sec:intro}

As many systems are too complex to be modelled (only) by physical relationships, control and system optimization procedures often need to employ data, making data-driven decision making a central topic. This topic can be considered as old as control itself. It is not hard to envisage that the flyball governor controlling Watt's steam engine was adjusted based on observations of the closed-loop operation of the engine. A more modern but still early example is the MIT rule, proposed for adaptive control of aircrafts traversing different flight conditions \cite{Whitaker&Yamron&Kezer:58}. Since then several prominent research directions have been pursued, focusing on different aspects of data-driven decision making in a control context. In adaptive control, different controller structures were proposed, such as the self-tuning regulator (STR) and model reference adaptive control (MRAC). These two principles employ data in a somewhat different way. In the STR a model is updated which subsequently is used to update the controller employing the certainty equivalence principle, i.e. the model is assumed to be correct, while in MRAC data directly affects the controller parameters. Establishing closed loop stability of adaptive control schemes (e.g. \cite{Egardt:79b,Goodwin&Ramadge&Caines:80}) is considered one of the breakthroughs in control. We refer to   \cite{Astrom&Wittenmark:95,Aström2014,Ioannou&Sun:96} for treatments of adaptive control. Robust control \cite{doyle89} put the spotlight on the approximative nature of models and the field "identification for control (I4C)" evolved in the 1990s as a response to the dichotomy between the leading paradigm of the time of using data-driven models as complex as the true system and the fact that simple controllers may successfully control a complex system. Control relevant models (Examples 2.10 and 2.11 in \cite{Astrom&Wittenmark:95} are striking illustrations) and procedures for how to generate data such that despite a systematic error (bias) the identified models fulfilled their task in control design were developed. We refer to \cite{Hjalmarsson:04a} for a survey of this work. Another line of research has been to develop systematic experiment design procedures such that the generated data is maximally informative for data-driven model based control design  \cite{Gevers&Ljung:86,Geversetal:06a,Hjalmarsson:09}. Interestingly, in \cite{Hjalmarsson:09} it is shown that even if such a design is made for a full order model, the experiments become control-relevant allowing for simplified models to be used as long as they can capture the system properties relevant for control. Other observations made in \cite{Hjalmarsson:09} are that the experimental cost and the required model complexity is highly dependent on the desired control performance. Also relevant to our study is that to cope with the inherent catch that the optimal experiment design depends on the unknown system, adaptive experiment design procedures have been developed supported by theory showing that such procedures do not result in loss in performance results asymptotically \cite{Huang:14b}. Control performance and the cost of acquiring data have in most work been treated separately. A seminal step for a comprehensive treatment of data-driven decision making for control was taken by Lai and Wei in the mid 1980s when they in an adaptive control context integrated exploration (purposely adding excitation to obtain informative data) and exploitation (achieving good control performance) into one single criterion \cite{Lai:86a,lai1986asymptotically}. The difference between the actual control performance, including both these contributions, and the optimal control performance was called regret. While dormant for long, lately significant developments have been made for the problem of regret minimization for the Linear Quadratic Regulator (LQR) problem. One of the main outcomes of these efforts has been to establish the rate of growth of the minimal regret as function of the control horizon $T$. The early work \cite{Lai:86a} indicated an asymptotic lower bound of $O(\log(T))$ for minimum variance control of ARX-systems. For LQR, the rate $O(\sqrt{T})$ was established in~\cite{ Mania2019CertaintyEI} and proven to be the optimal lower bound rate when both state matrix and input matrix are unknown in~\cite{pmlr-v119-simchowitz20a}. 
This rate can be obtained when systems are excited with white Gaussian noise whose variance decays as $1/\sqrt{t}$~\cite{Wan:21}. We will refer to this type of exploration as decaying. Recently, in~\cite{jedra2021minimal}, it is demonstrated that when both state matrix and input matrix are unknown, the regret can be upper-bounded as $O(\sqrt{T})$;  when either state matrix or input matrix is known, it can be upper-bounded as $O(\log(T))$. The exploration cost is also accounted for in  \cite{Forgione:15}, who studies a general control problem for a general class of linear systems in a batch-wise setting. The numerical study shows that it is optimal to devote all exploration to the first batch. Even though this setting is different from the adaptive LQR setting studied in the works referenced above, it is interesting to note that this conclusion is at odds with the decaying exploration proposed in, e.g. \cite{Wan:21,jedra2021minimal} and in~\cite{colin2023optimal} it is shown that for a given finite horizon $T$ immediate excitation is optimal also in an LQR-setting. By immediate we mean that all exploration takes place at the beginning. 
The optimal regret still is $O(\sqrt{T})$ but the proportionality constant is smaller than with decaying. 

In this contribution we broaden this line of research by turning to the optimization of non-linear static systems where the performance measure can be a non-linear function of the input. To study essential aspects of the problem \textcolor{black}{and address the conceptual differences and challenges that arise in the non-linear context}, we limit our attention to static input-output relationships described by one unknown parameter and assume that we have access to noisy measurements of the output response. This does not preclude that the system may be dynamic, only that we are interested in optimizing its static behaviour and that we assume that the response to a constant input can be measured (for a stable system the latter may be achieved after transients have died out). The setting thus closely relates to the steady-state optimization problems faced in Real-Time Optimization (RTO) \cite{Marchetti2016, Srinivasan2019} where the goal is to optimize a plant operating condition, despite such plant's input-output steady-state map being partially unknown, e.g. as it depends on an unknown parameter to be learned with noisy measurements. Although the concept of exploration-exploitation trade-off has been central in the context of decision making under uncertainty, only few works in the literature of RTO account for this. A notable exception is  \cite{delRio2021} which employs a Bayesian Optimization framework, so that exploration is accounted for by means of an acquisition function. To the best of our knowledge, a regret minimization framework was, for this type of setting, first studied in \cite{pasquini2023e2}. However, 
there the long-term effect of the additional excitation is neglected as only the regret in the next time instant is considered. Contrary to this, we derive an approximation for the cumulative regret over the entire time horizon $T$, based on which we are able to derive an \textcolor{black}{effective exploration} strategy.
Aligning with the numerical results of \cite{Forgione:15} for batch-wise regret minimization 
and the theoretical results of \cite{colin:23a,colin2023optimal} for the LQR-problem, we show that \textcolor{black}{an effective exploration} strategy is either lazy (no exploration at all) or immediate (only explore at the first time instant), depending on the problem at hand. 


%

\vspace{-0.1cm}\section{Problem statement}
\label{sec:method}

Many important control problems can be phrased as iterative optimization of the input. One example is medium optimization in bioprocessing applications for pharmaceuticals production \cite{Wang2024}. More broadly, this is the scope of real-time optimization, where the operating condition of a partially unknown plant is optimized by iteratively solving a sequence of optimization problems. Here we consider the problem of unconstrained optimization of a scalar cost function $\Phi$ where the dependency on the underlying system is modeled by a scalar parameter $\theta_0$, i.e. $\Phi=\Phi(u,\theta_0)$ where $u$ denotes the scalar input to be optimized and $\theta_0$ denotes the true parameter. The optimal input is given by 
\begin{equation}\label{eq:exploit}
	u_0^*=\arg\min_{u\in\mathbb{R}} \Phi(u, \theta_0)
\end{equation} 
We use the function $U:\mathbb{R}\rightarrow \mathbb{R}$ to explicitly describe the relationship between the system parameter and the optimal input defined by~\eqref{eq:exploit}, i.e. $u_0^*= U(\theta_0)$.
A key aspect of our problem is that exact knowledge about $\theta_0$ is unknown  and is only available indirectly by way of measurements of some output of the system subject to measurement errors. Formally, we express this by the measurement equation
\begin{equation}\label{eq:noisy}
	y_t = h(u_t, \theta_0)+e_t
\end{equation}
where $y_t$ and $u_t$ are the output and input at time instant $t$ respectively
and $e_t$ is zero-mean Gaussian noise with variance $\sigma^2$.
We will base our method on the Certainty Equivalence Principle (CEP), i.e. after having collected $\{u_1,\dots,u_{t-1}\}$ and $\{y_1,\dots,y_{t-1}\}$ from \eqref{eq:noisy},  we approximate $u_0^*$ by replacing the unknown $\theta_0$ by an estimate $\hat{\theta}_t$ obtained using those measurements, resulting in the input
\begin{equation}\label{eq:exploit2}
	u_t^*=\arg\min_{u} \Phi(u, \hat{\theta}_t)
\end{equation}
We remark that the use of CEP is very common in data-driven control \cite{Astrom&Wittenmark:95,Wan:21,Mania2019CertaintyEI}.
  
From~\eqref{eq:exploit} and~\eqref{eq:exploit2} we see that $u_t^*=U(\hat{\theta}_t)$. 
Thus, attaining an accurate estimate is  instrumental for the CEP to give a satisfactory result. To this end the CEP may not work well as the input may not be very informative. 
Consider for example the extreme case for which $u=u_0^*$ gives $y=0$, i.e. no information on $\theta_0$ is obtained. 

\vspace{-0.4cm}
\begin{figure}[h!]
        \begin{center}		           \includegraphics[width=0.8\columnwidth]{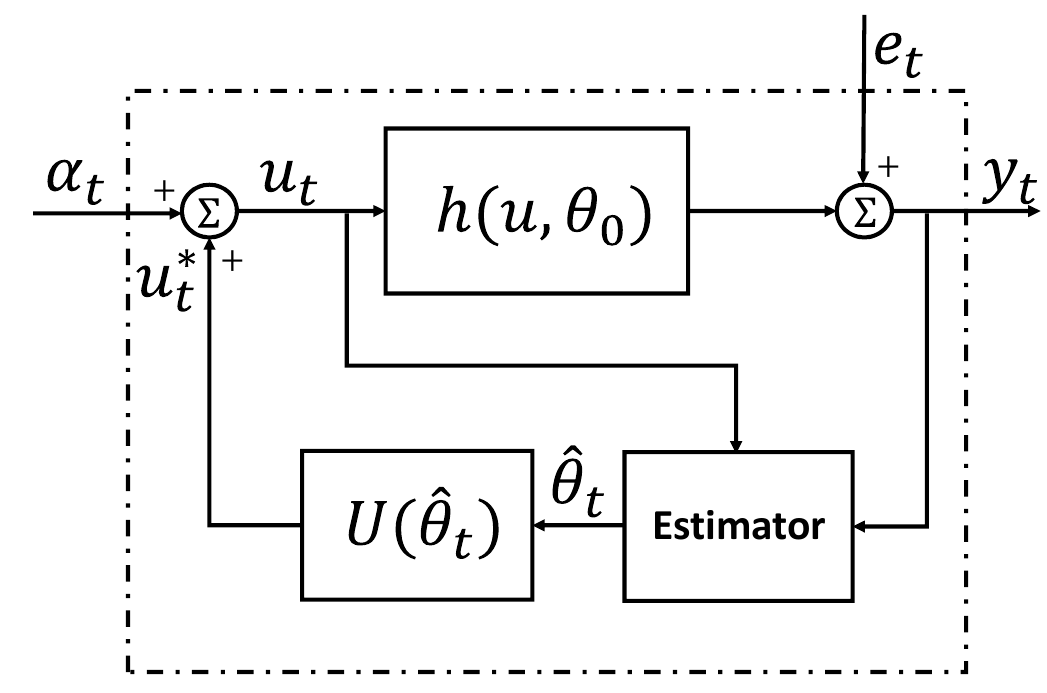}    
        \caption{The iterative framework for the input optimization problem based on the CEP as well as the exploitation and exploration idea} 
        \label{fig:diagram}
	\end{center}
 \vspace{-0.6cm}
 
\end{figure}

To address this, following e.g. \cite{Wan:21}, exploration is achieved by incorporating an additional excitation term $\alpha_t$, to the input $u_t^*$ to improve the accuracy of estimation, as illustrated in Fig.~\ref{fig:diagram}. For future reference we will refer to $u_t^*$ and $\alpha_t$  as exploitation and exploration inputs respectively. However, the introduction of $\alpha_t$ presents a dilemma, as it tends to perturb $u_t^*$, the input that minimizes $\Phi(u, \hat{\theta}_t)$. This perturbation leads to an increase in the cost. Thus, a trade-off exists between exploitation and exploration when we design $\alpha_t$.

Regret is a widely used criterion for obtaining an optimal balance between exploitation and exploration.
It is defined as the cumulative performance degradation from the optimal cost achieved with $u_0^*$, when the input $u_t = u^*_t + \alpha_t$ is employed instead. 
In a finite horizon $T$, the cumulative regret is expressed as $\sum_{t=1}^{T} \left[ \Phi(u^*_t + \alpha_t, \theta_0) - \Phi(u^*_0, \theta_0) \right]$ which can be used to guide the design of the exploration $\alpha_t$.

The presence of random noise $e_t$, as described in equation~\eqref{eq:noisy}, 
makes both $\hat{\theta}_t$ and the subsequently generated $u^*_t$ stochastic.
Inspired by the literature on regret minimization for linear quadratic adaptive control~\cite{Wan:21,jedra2021minimal,colin2023optimal}, \textcolor{black}{we choose the exploration sequence $\{\alpha_t\}$ to be zero-mean white noise with a time-varying variance denoted by $x_t$.}
Thus, our further analysis will focus on the expected cumulative regret (in what follows we will use regret for short for this quantity) defined as follows, taking into consideration the randomness in both the noise and the exploration action,
\begin{equation}\label{eq:regret_alpha}
	\bar{R}=\sum_{t=1}^{T}
	\mathbb{E}\big[
	\Phi(u^*_t+\alpha_t, \theta_0)-\Phi(u^*_0, \theta_0)
	\big]
\end{equation}
Here the expectation $\mathbb{E}$ is taken with respect to (w.r.t.) both the noise $\{e_t\}$ and the exploration signal $\{\alpha_t\}$. This reflects the average performance loss of the exploration strategy when the same experiment is repeated multiple times with the same way of generating the stochastic $e_t$ and $\alpha_t$.


In general it is very difficult to compute the regret in \eqref{eq:regret_alpha}. In the next section we will develop an approximate expression of $\bar{R}$ which will allow us to conclude on the structure of the \textcolor{black}{optimal exploration} sequence $\{\alpha_t\}$.

\section{Regret approximation \& model uncertainties}\label{sec:approx}

\vspace{-0.08cm}\subsection{Regret approximation}\label{subsec_reg}
\vspace{-0.081cm}

\textcolor{black}{We assume that the exploitation input $u_t^*$ is in the vicinity of $u_0^*$ and hence we use a  Taylor expansion of~\eqref{eq:regret_alpha} around $u_0^*$}
\begin{align}\label{eq:original R} 
	\bar{R}\approx \sum_{t=1}^{T}
	\mathbb{E}\big[	
	\big(u^*_t&+\alpha_t-u^*_0\big)J_u
	+\frac{H_u}{2}
	\big(u^*_t+\alpha_t-u^*_0\big)^2
	\big]
\end{align}
where $J_u=\frac{\partial \Phi(u, \theta_0)}{\partial u}|_{u=u^*_0}$ and $ 
	H_u=\frac{\partial^2 \Phi(u, \theta_0)}{\partial u^2}|_{u=u^*_0}$.
Here, $J_u=0$ since $u^*_0$ is a minimum of $\Phi(\cdot, \theta_0)$. For simplicity we will assume that $H_u>0$. Since  $H_u/2$ is just a scaling factor we will omit it and study the scaled expected regret $\sum_{t=1}^{T}\mathbb{E}
	\big[(u^*_t+\alpha_t-u^*_0)^2\big]$, which can be expanded into
\begin{align}\label{eq:taylor R}
        \sum_{t=1}^{T}
	\mathbb{E}\big[
	(u_t^*-u^*_0)^2+2(u_t^*-u^*_0)\alpha_t+\alpha^2_t
	\big]
\end{align}
Furthermore, both the actual input $u^*_0$ and the approximate counterpart $u_t^*$ are determined by the function $U$.
The former requires the unknown actual parameter $\theta_0$, while the latter utilizes the estimated parameter $\hat{\theta}_{t}$.
We here approximate the difference between $u^*_0$ and $u_t^*$ with the Taylor approximation
\vspace{-0.9cm}

\begin{align}\label{eq:taylor2}
	u_t^*-u^*_0\approx 
	(\hat{\theta}_{t}-\theta_0)J_{\theta}
\end{align}
where $J_{\theta}=\frac{dU(\theta)}{d\theta}|_{\theta=\theta_0}$. 
This gives the following approximation of the regret in~\eqref{eq:taylor R} 
\begin{align}\label{eq:exp_R}
	\tilde{R}:= &
	\sum_{t=1}^{T}
	\mathbb{E}\big[
	(\hat{\theta}_{t}-\theta_0)^2J^2_{\theta}+2(\hat{\theta}_{t}-\theta_0)J_{\theta}\alpha_t+\alpha^2_t
	\big] \nonumber \\	=&J^2_{\theta}\sum_{t=1}^{T}\mathbb{E}\big[
	(\hat{\theta}_{t}-\theta_0)^2\big]
	+\sum_{t=1}^{T}x_t 
\end{align}
where the equality holds because the current exploration action $\alpha_t$ and the estimate error $\hat{\theta}_{t}-\theta_0$ (determined by the input and output before time $t$) are independent, together with the fact that $\mathbb{E}[\alpha_t]=0$. 
The approximation \eqref{eq:exp_R} indicates that the performance degradation over horizon $T$ is incurred by the model uncertainty and the exploration input in an additive way. Next we will study the model uncertainty and its relation to the exploration input. 


\vspace{-0.08cm}\subsection{Model uncertainty}
\label{subsec:fi}
\vspace{-0.081cm}

The Cramér-Rao inequality establishes a lower bound for the variance of an unbiased estimator in terms of the inverse of the Fisher information~\cite{Ljung:1999}. \textcolor{black}{The lower bound is reached\footnote{This holds when $h$ is linear in $\theta$ under no feedback but is otherwise an idealization. For large $t$ under sufficient excitation it holds in general.} when we adopt the following idealized assumption 
\begin{assumption}\label{as0}
\textcolor{black}{The estimator $\hat{\theta}_{t}, \ t=1, \cdots, T$, is unbiased and efficient~\cite{Ljung:1999}.}
\end{assumption}}
\noindent Thus, the estimate variance $\mathbb{E}\big[(\hat{\theta}_{t}-\theta_0)^2\big]$ in~\eqref{eq:exp_R} is given by the inverse of the Fisher information denoted by $\mathbb{I}_{t-1}$, where the index is $t-1$ since the model uncertainty will depend on all the inputs up to, but not including, time $t$.
\textcolor{black}{As shown in Appendix~\ref{sec:app1}, the Fisher information at time $t$ is expressed as follows}
\begin{align}\label{eq:Fish_t}
	\mathbb{I}_{t} 
	=& \mathbb{I}_{0} + 
        \frac{1}
	{\sigma^{2}}
	\sum_{s=1}^{t}	\mathbb{E}\bigg[	
	\frac{\partial h}{\partial \theta}	 \Bigl|_{\substack{\theta=\theta_0\\u_s=u^*_s+\alpha_s}}^2
	\bigg] 
\end{align}
where the initial information, denoted as $\mathbb{I}_0$, is obtained from prior experiments.
In~\eqref{eq:Fish_t}, the input $u_s$ consists of two parts, the exploration input $\alpha_s$ and the exploitation input $u^*_s$ as illustrated in Fig.~\ref{fig:diagram}. 
We here introduce a \textcolor{black}{simplified toy example that is useful to continue our discussion and illustrate the new findings}. \textcolor{black}{This is intended not as a comprehensive analysis, but as a clear demonstration of our ideas, and it} will be further developed in Section \ref{sec:example}.

\begin{example}\label{example}
    Let us consider the following quadratic cost function and input-output relationship 
\begin{equation}\label{eq:eg}
	\begin{aligned}
		&\Phi(u, \theta_0) = u^2 + 2(\theta_0+1)u \\
		&y = h(u, \theta_0)+e=\theta_0u^2+e
	\end{aligned}
\end{equation}
where $e$ is zero-mean white Gaussian noise with $\sigma^2=1$. The optimal, but unknown, input $u_0^*$ for the quadratic function $\Phi$ is given by
$u_0^*= U(\theta_0)=-(\theta_0 + 1)$. We consider the iterative framework presented in Fig \ref{fig:diagram}, where the input $u_t$ applied at iteration $t$ is composed of an exploitation input $u_t^* = -(\hat\theta_t + 1)$ (via the CEP) and an exploration input $\alpha_t$. With the notation of Section~\ref{subsec_reg}, we notice that in this case  $J_\theta$ is constant since $J_{\theta}=\frac{dU(\theta)}{d\theta}|_{\theta=\theta_0} = -1$. Then, under Assumption \ref{as0} that the estimator is unbiased and efficient, following \eqref{eq:exp_R} and \eqref{eq:Fish_t}, we have
\begin{align}\label{eq:final_eq}
	\tilde{R}=& \sum_{t=1}^{T}
	\frac{1}{\mathbb{I}_{t-1}}
	+
	\sum_{t=1}^{T}x_t = 
        \frac{1}{\mathbb{I}_0}
        +
	\sum_{t=1}^{T-1}
	\frac{1}{\mathbb{I}_t}
	+
	\sum_{t=1}^{T}x_t \\
	=&
        \frac{1}{\mathbb{I}_0}
        +
	\sum_{t=1}^{T-1}
	\frac{1}{
		\mathbb{I}_0+
		\sum_{s=1}^{t}
		\mathbb{E}
		\big[
		u_s^4|_{\substack{u_s=u^*_s +\alpha_s}}
		\big]
	}
	+
	\sum_{t=1}^{T}
	x_t \nonumber
\end{align}

For the white noise exploration input $\alpha_t$, let us choose  a zero-mean Gaussian distribution with time-varying variance $x_t$. The sequence $\{x_t\}$ contains the decision variables to be designed for regret minimization.
We can now expand 
\begin{align}\label{eq:fi_eg}
    & \mathbb{E}
	\big[(u_s^*+\alpha_s)^4
	\big] 
	=
	\mathbb{E}
	\big[\big(u^*_0+(\hat{\theta}_s-\theta_0)J_{\theta}+\alpha_s\big)^4
	\big]  
\end{align}where the equality comes from the Taylor approximation in~\eqref{eq:taylor2} which in this example is exact since $u_s^* = -(\hat\theta_s + 1)$ is linear w.r.t. $\hat\theta_s$. By expanding the right hand side of the latter, using the independence between $\alpha_s$ and the estimate error $\hat{\theta}_s-\theta_0$ (determined by the input before time $s$), and recalling that $\alpha_s$ is zero-mean Gaussian with variance $x_s$, and that $\hat{\theta}_s - \theta_0$ is zero-mean Gaussian (this holds when there is no feedback), we get 
\begin{align}\label{eq:fi_eg2}
    & 3x_s^2+\big[6\mathbb{E}[(\hat{\theta}_s-\theta_0)^2]+6(u^*_0)^2\big]x_s+\mathbb{E}[(\hat{\theta}_s-\theta_0)^4]\\ 
     &+6(u^*_0)^2\mathbb{E}[(\hat{\theta}_s-\theta_0)^2]+(u^*_0)^4  
\end{align}Finally, by recalling that $\hat{\theta}_s$ is assumed to be a unbiased and efficient estimator for $\theta_0$ which implies that the variance of $\hat{\theta}_s - \theta_0$ is equal to $1/\mathbb{I}_{s-1}$ and every odd moment of $\hat{\theta}_s-\theta_0$ is zero, we get the final expression for $  \mathbb{E}
	\big[(u_s^*+\alpha_s)^4
	\big] $
\begin{equation}
\label{eq:last_expansion_example_information}
  3x_s^2+[6\mathbb{I}^{-1}_{s-1}+6(u^*_0)^2]x_s+3\mathbb{I}^{-2}_{s-1}+6(u^*_0)^2\mathbb{I}^{-1}_{s-1}+(u^*_0)^4 
\end{equation}
We can then rewrite \eqref{eq:final_eq} as 
\begin{align}\label{eq:final_eq_with_I_cal}
	\tilde{R}=&
        \frac{1}{\mathbb{I}_0}
        +
	\sum_{t=1}^{T-1}
	\frac{1}{
		\mathbb{I}_0+
		\sum_{s=1}^{t}
		\mathcal{I}(x_s,\mathbb{I}_{s-1}^{-1})
	}
	+
	\sum_{t=1}^{T}
	x_t  
\end{align}
where $\mathcal{I}(x_s,\mathbb{I}_{s-1}^{-1})$ is the expression \eqref{eq:last_expansion_example_information}, which is a function of the exploration variance at time instant $s$ and the inverse of the Fisher information at time instant $s-1$.
\end{example}

As a further simplification we will approximate the terms in \eqref{eq:Fish_t} using the approximation of $u^*_s$ as introduced in~\eqref{eq:taylor2} to get the approximate Fisher information
\begin{align}\label{eq:app_Fish_t}
	\tilde{\mathbb{I}}_{t} 
	=& \mathbb{I}_{0} + 
        \frac{1}
	{\sigma^{2}}
	\sum_{s=1}^{t-1}	\mathbb{E}\bigg[	
	\frac{\partial h}{\partial \theta}	        \Bigl|_{\substack{\theta=\theta_0\\u_s= u^*_0+(\hat{\theta}_s-\theta_0)J_{\theta}+\alpha_s}}^2
	\bigg] 
\end{align}
As we saw in Example \ref{example}, each term in the sum of the approximate information \eqref{eq:app_Fish_t} (which was exact in Example~\ref{example}) is dependent on the exploration variance $x_s$ and the inverse of the Fisher information $\mathbb{I}_{s-1}$. This justifies the formal introduction of the incremental information function  $\mathcal{I}:\mathbb{R}_+^2\rightarrow \mathbb{R}_+$, defined as 
\begin{align}\label{eq:info_fun}
	\mathcal{I}(x_s, \mathbb{I}_{s-1}^{-1}) \coloneqq
        \frac{1}{\sigma^{2}}\mathbb{E}\bigg[    \frac{\partial h}{\partial \theta}        \Big|_{\substack{\theta=\theta_0\\u_s= u^*_0+(\hat{\theta}_s-\theta_0)J_{\theta}+\alpha_s}}^2  \bigg]
\end{align}where $\mathbb{R}_+$ is the set of non-negative real scalars. 
\begin{remark}
Example \ref{example} presented the case of an incremental information function $\mathcal{I}$ dependent on the exploration variance $x_s$ and the inverse of Fisher information at time instant $s-1$. We argue that this is the case in general, as the square of the partial derivative in~\eqref{eq:info_fun} can  either be directly expanded, or approximated with arbitrary precision with Taylor expansion, to be a polynomial function, consisting of terms proportional to $\alpha_s^n(\hat{\theta}_s-\theta_0)^m$, where $n$ and $m$ are all possible exponents combinations.
When taking the expectation in \eqref{eq:info_fun}, the independence between the action $\alpha_s$ and the estimate error $\hat{\theta}_s-\theta_0$ (determined by the inputs before time $s$), allows to further simplify this expression by using $\mathbb{E}[\alpha_s^n(\hat{\theta}_s-\theta_0)^m]=\mathbb{E}[\alpha_s^n]\mathbb{E}[(\hat{\theta}_s-\theta_0)^m]$. This, together with \textcolor{black}{Assumption \ref{as0} of an} unbiased and efficient estimator, justifies the dependency of the incremental information function on the action variances $x_s$ and the inverse of the Fisher information at time $s-1$, i.e. $\mathbb{I}_{s-1}$.
\end{remark}
From \eqref{eq:info_fun} we can rewrite the approximate Fisher information in \eqref{eq:app_Fish_t} as
\begin{align}\label{eq:I_at_t}
	\tilde{\mathbb{I}}_{t} 
	= \mathbb{I}_{0} + \sum_{s=1}^{t}	 \mathcal{I}(x_s, \mathbb{I}^{-1}_{s-1})=\tilde{\mathbb{I}}_{t-1}+ \mathcal{I}(x_t, \mathbb{I}^{-1}_{t-1})
\end{align}
where $\tilde{\mathbb{I}}_{t-1}={\mathbb I}_0 + \sum_{s=1}^{t-1}	 \mathcal{I}(x_s, \mathbb{I}^{-1}_{s-1})$, which in turn gives the following approximation of the regret~\eqref{eq:exp_R} 
\begin{align}
	\tilde{R} =
 &\sum_{t=1}^{T}\frac{J^2_{\theta}}{\mathbb{I}_{t-1}}
	+\sum_{t=1}^{T}x_t 
 \approx 
 \frac{J^2_{\theta}}{\mathbb{I}_{0}}+\sum_{t=1}^{T-1}\frac{J^2_{\theta}}{\tilde{\mathbb{I}}_{t}}
	+\sum_{t=1}^{T}x_t 
 \nonumber
 \\\label{eq:bound R} = &\frac{J^2_{\theta}}{\mathbb{I}_{0}}+\sum_{t=1}^{T-1}\frac{J^2_{\theta}}{\mathbb{I}_{0} + \sum_{s=1}^{t}	 \mathcal{I}(x_s, \mathbb{I}^{-1}_{s-1}) }
	+\sum_{t=1}^{T}x_t 
\end{align}	
The non-linear dynamics of the approximate Fisher information in \eqref{eq:I_at_t} presents the main challenge in devising an \textcolor{black}{optimal exploration} strategy with minimizing the regret. 
The incremental information function $\mathcal{I}$ varies with the cost function $\Phi(u, \theta_0)$ and the input-output relationship $h(u, \theta_0)$. In the following we will focus on the case where the incremental information function $\mathcal{I}$ satisfies the following assumption (which is verified, e.g., in Example \ref{example}).
\begin{assumption}\label{as1}
 $\mathcal{I}$ is non-negative, monotonically increasing and convex w.r.t. the first argument.
 Furthermore, $\mathcal{I}$ is monotonically increasing w.r.t. the second argument.
\end{assumption}
Under Assumption \ref{as1} it holds that
\begin{align}	
	\tilde{\mathbb{I}}_{t}= \mathbb{I}_{0} + \sum_{s=1}^{t}	 \mathcal{I}(x_s, \mathbb{I}^{-1}_{s-1})   > \mathbb{I}_{0} +  \sum_{s=1}^{t}	 \mathcal{I}(x_s, 0)
\end{align}
Using this, the approximate regret in~\eqref{eq:bound R} can be  upper-bounded by
\begin{equation}\label{eq:regret_ub}
	\begin{aligned}
		R_{ub} \coloneqq \frac{J^2_{\theta}}{\mathbb{I}_{0}}+\sum_{t=1}^{T-1}\frac{J^2_{\theta}}{\mathbb{I}_{0} +\sum_{s=1}^{t} \mathcal{I}(x_s, 0)}+\sum_{t=1}^{T}x_t
	\end{aligned}
\end{equation}
In the following, we use the upper bound of the approximate regret as a design criterion for developing an excitation strategy.
Before this we will rewrite the expression for $R_{ub}$ to simplify the notation slightly. 
Given that $\mathcal{I}$ in~\eqref{eq:regret_ub} has a constant second argument, we define a new information function $i:\mathbb{R}\rightarrow \mathbb{R}$ as $i(x_s) \coloneqq \mathcal{I}(x_s,0)/J_{\theta}^2$, which inherits the properties of $\mathcal{I}$ w.r.t. its first argument, i.e. non-negative, monotonically increasing and convex.
Besides, we notice that $x_T$ only appears in the last sum in \eqref{eq:regret_ub} and therefore it is optimal to take $x_T=0$.
By defining $i_0\coloneqq\mathbb{I}_{0}/J_{\theta}^2$, the upper bound~\eqref{eq:regret_ub} is rewritten as
\begin{equation}\label{eq:final R}
	\begin{aligned}
		R_{ub}=\frac{1}{i_0} + \sum_{t=1}^{T-1}
	\frac{1}{
	i_0+\sum_{s=1}^{t}
	i(x_s)}
+\sum_{t=1}^{T-1}x_t
	\end{aligned}
\end{equation}

\vspace{-0.1cm}\section{Theoretical result}
\label{sec:theoretical_result}

Our next task is to minimize \eqref{eq:final R} w.r.t. 
the vector $x = [x_1, \dots, x_{T-1}]$, which consists of the non-negative variances of the exploration input at each time step. 
We notice that the first term in \eqref{eq:final R} is constant and can be omitted.
Before stating a theorem for this problem, we introduce two families of excitation signals.

\begin{definition}
	\label{def1}
	We say that $x \in \mathbb{R}^{T-1}$ is an immediate excitation if 
		$x_1>0 \text{ and }  x_2=\cdots=x_{T-1} = 0$, and a lazy excitation if
		$x_k = 0, \ k = 1, \dots, T-1$.

	
\end{definition}

\begin{theorem}
	\label{th1}
	Consider the problem
        \begin{align}         \label{eq:optimization_problem_theorem}   
            \min\limits_{x_1,\dots,x_{T-1}} \ \ &\sum_{t=1}^{T-1}
	\frac{1}{
	i_0+\sum_{s=1}^{t}
	i(x_s)}
+\sum_{t=1}^{T-1}x_t \\ \rm{s.t.} \quad \ & x_k \ge 0, \quad k = 1,\dots,T-1 \nonumber
        \end{align}
 and assume that the information function $i(\cdot)$  is non-negative, monotonically increasing and convex in the domain $[0, \infty)$. 
	Let $x^*$ be the optimal solution of \eqref{eq:optimization_problem_theorem}. 
	Then $x^*$ is either a lazy or an immediate excitation (see Definition~\ref{def1}).	
		Moreover, if the following inequality holds 
\begin{equation}\label{eq:immed_condition}
	\sum_{t=1}^{T-1}\frac{i'(0)}{[i_0+t\cdot i(0)]^2}> 1.
\end{equation}
then $x^*$ is an immediate exploration solution.
\end{theorem}

\noindent \textit{Proof}. \textcolor{black}{See Appendix~\ref{sec:app2}.}

\begin{remark}\label{rem:interpretation_th1}
The intuitive interpretation of Theorem~\ref{th1}  is that when exploration is necessary, it is best to do it as early as possible since the reward in terms of lower cost, due to a better model, then accumulates over the entire horizon $T$, rather than a portion of it. The reduction in cost that can be achieved is also higher early on since the information in data then is less than later on (recall that the information at a certain time depends on data up to that time point). 

Moreover, the findings of Theorem~\ref{th1} \textcolor{black}{for a static, scalar and non-linear problem} strongly resonate with the numerical results of~\cite{Forgione:15} and the theory  in~\cite{colin2023optimal,colin:23a} for the LQR problem where  the exploration strategy is either a lazy or an immediate excitation. Theorem~\ref{th1} goes beyond these works since the Fisher information is not restricted to be linear w.r.t. the exploration decision variable as is the case in~\cite{Forgione:15,colin2023optimal,colin:23a}.
\end{remark}
\begin{remark}\label{rem:condition} \textcolor{black}{
The sufficient condition in~\eqref{eq:immed_condition} for immediate excitation suggests that immediate exploration is more likely under the following cases: (1) a large horizon $T$; (2) a large value of $i'(0)$, implying that even a small exploration yields significant new information; (3) a small initial information $i_0$; (4) a small information $i(0)$  in the absence of  exploration.
} \vspace{-0.3cm} \end{remark}

\vspace{-0.1cm}\section{Numerical example (continued)}
\label{sec:example}

\vspace{-0.08cm}\subsection{Objectives of the example}
\vspace{-0.081cm}

In this section, we return to Example~\ref{th1}. The purpose is to check if a design of the exploration based on $R_{ub}$ in~\eqref{eq:final R} (by taking advantage of the results of Theorem~\ref{th1})  leads to satisfactory performance for the actual regret $\bar{R}$ in~\eqref{eq:regret_alpha}. We will consider explorations of the form:
\begin{equation}\label{eq:white_noise_excitation}
    \alpha_t = \sqrt{x_t} \bar{\alpha}_t
\end{equation}
where $\{\bar{\alpha}_t\}$ is a zero-mean white noise and $\{x_t\}$ is the variance sequence to be designed. We consider the following four particular cases of such kind of excitation 
\begin{itemize}
    \item immediate Gaussian: each $\bar{\alpha}_t$ is drawn from a zero-mean normal distribution with unit variance, $x_2 = \cdots = x_T = 0$  and $x_1 = x_g$ where $x_g \ge 0$ is to be tuned.
    \item immediate binary: each $\bar{\alpha}_t$ is drawn from a zero-mean binary distribution whose two possible values are $-1$ and 1, $x_2 = \cdots = x_T = 0$  and $x_1= x_b$ where $x_b \ge 0$ is to be tuned.
    \item lazy exploration: $x_1 = \cdots = x_T = 0$.
    \item decaying Gaussian: each $\bar{\alpha}_t$ is drawn from a zero-mean normal distribution with unit variance and  $x_t$ decays as $  x_t = c t^p$
where $c\ge 0$ and $p<0$.  This choice \textcolor{black}{with~$p=-0.5$} is frequently used in the LQR literature~\cite{Wan:21,jedra2021minimal}.
\end{itemize}
To validate the effectiveness of the \textcolor{black}{exploration} strategy in Theorem~\ref{th1} (based on the minimization of $R_{ub}$)   for the minimization of the actual regret $\bar{R}$, we compare the actual regret $\bar{R}$ obtained with these four exploration strategies with well-tuned values $x_g$, $x_b$, $c$ and $p$ obtained by:
\begin{itemize}
  \item[$(a)$]  minimizing $R_{ub}$.
  \item[$(b)$]  directly minimizing the actual regret $\bar{R}$. 
\end{itemize}The purpose of design $(b)$ is to illustrate the error  incurred by the exploration design with minimizing ${R}_{ub}$ instead of with $\bar{R}$. The term $R_{ub}$ has the form~\eqref{eq:final R} where the information function $i(x_s)$ for this example is given by:
\begin{itemize}
    \item $i(x_s)=3x_s^2+6(u^*_0)^2x_s+(u^*_0)^4$ when $\bar{\alpha}_t$ in~\eqref{eq:white_noise_excitation} is drawn from a zero-mean Gaussian distribution.
    \item $i(x_s)=x_s^2+6(u^*_0)^2x_s+(u^*_0)^4$ when $\bar{\alpha}_t$ in~\eqref{eq:white_noise_excitation} is drawn from a zero-mean binary distribution.
\end{itemize} In both cases the information function $i(x_s)$ is a non-negative, monotonically increasing and convex function in $[0,+\infty)$. Hence, according to Theorem~\ref{th1}, the exploration vector 
 $x^* = [x_1^*, \dots, x_{T-1}^*]$ minimizing $R_{ub}$ is either lazy or immediate for both distribution choices for $\{\bar{\alpha}_t\}$. 

\vspace{-0.08cm}\subsection{Simulation details}
\vspace{-0.081cm}

We will consider a horizon of $T = 50$ and choose $\sigma^2 = 1$ assumed to be known\footnote{It can be estimated together with $\hat{\theta}_t$ if not known.}. In order to implement the proposed scheme, we require an initial estimate for $\theta_0$ and the initial Fisher information $\mathbb{I}_0$. For this purpose, we perform an initial identification with one input-output data pair collected on the system excited with a deterministic input equal to $1$. 

We conduct the simulation on 10 different systems, each characterized by a different parameter $\theta_0$ with the following values $\{-2,-0.7,-0.5,-0.4,-0.3,0.2,0.4,0.7,1,3\}$. 

For both designs $(a)$ and $(b)$, we will consider a grid-search approach in order to search for the optimal $x_g$, $x_b$, $c$ and $p$  with the following grid specifications:
\begin{itemize}
    \item Constants $x_g$, $x_b$ and $c$: 301 points log-regularly spaced  between  $10^{-3}$ and $10^{2}$.
      \item Exponent $p$: 21  points log-regularly spaced  between  $-20$ and $-0.1$.
\end{itemize}

  We will consider $N_{mc} = 1000$ Monte Carlo simulations to approximate $\Bar{R}$ in~\eqref{eq:regret_alpha}, which involves an expectation w.r.t. both $\{e_t\}$ and $\{\alpha_t\}$.   To ensure a fair comparison between all the exploration strategies,  we use \textit{the same} $N_{mc}$ realizations of a zero-mean white Gaussian noise with unit variance for ${e}_t$. Similarly, for the two distribution choices of $\{\bar{\alpha}_t\}$, we use \textit{the same} $N_{mc}$ zero-mean white Gaussian noise realizations with unit variance and \textit{the same} $N_{mc}$ binary sequences with unit variance. 
  
  For each value of $\theta_0$ and for each exploration strategy across all possible values of $x_g$, $x_b$, $c$ and $p$, we compute the average of the regret $\sum_{t=1}^T(\Phi(u^*_t+\alpha_t, \theta_0)-\Phi(u^*_0, \theta_0))$ obtained with the different realizations and the corresponding $R_{ub}$.  Then, we select the optimal parameters $x_g$, $x_b$, $c$ and $p$ that minimize $R_{ub}$ in design $(a)$ and $\bar{R}$ in design $(b)$.

\vspace{-0.08cm}\subsection{Results on the 10 systems}
\vspace{-0.081cm}

For the 10 systems,  we observe the following:
\begin{itemize}
    \item The lazy exploration never minimized $\bar{R}$. This is in line with our expectation, since for all $\theta_0$ values,  the sufficient condition~\eqref{eq:immed_condition}  for the optimality of immediate excitation in Theorem~\ref{th1}  was satisfied.
    \item Immediate binary provided the optimal regret $\bar{R}$ regardless of the design $(a)$ or $(b)$.
\end{itemize}For 8 out of the 10 systems, decaying Gaussian explorations resulted in a lower regret $\bar{R}$ than immediate Gaussian explorations, for both designs $(a)$ and $(b)$. This seems to contradict our theory. However, the exponent $p$ chosen for these 8 cases was -2.402, causing a rapid decrease in the variance $x_t$, such that the optimal decaying Gaussian explorations resemble the immediate ones. \textcolor{black}{Taking the excitation profiles for the system with $\theta_0 = 0.2$ as an example, we observe that the exploration variance for the decaying Gaussian exploration diminishes rapidly, closely resembling that of the immediate Gaussian exploration. The two strategies achieve the similar expected regret, with the decaying Gaussian performing slightly better. This observed difference can be attributed to approximations made during the analysis.}

\begin{figure}
    \centering
    \includegraphics[width =0.95\linewidth]{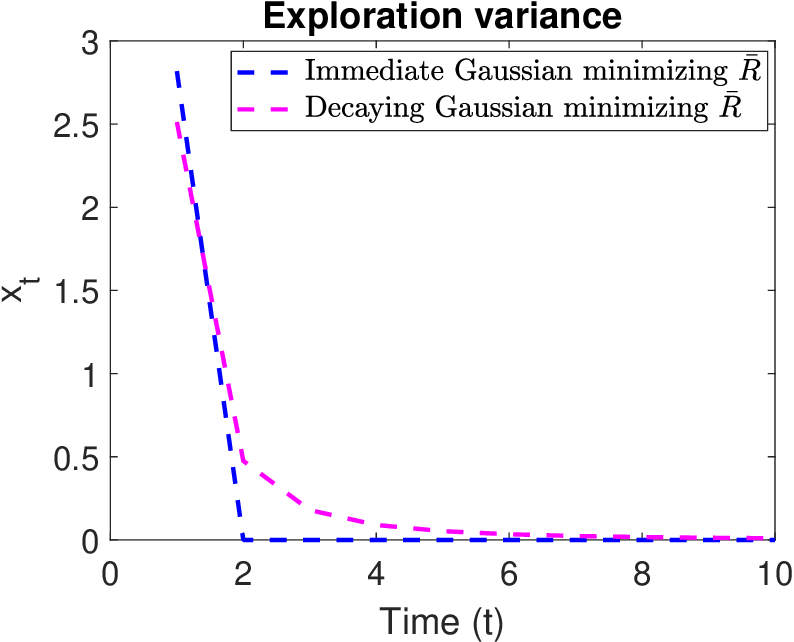}
    \vspace{-0.3cm}
    \caption{\textcolor{black}{Time evolution of the exploration variance with immediate Gaussian (blue line) and decaying Gaussian (magenta line), tuned with design $(b)$, for the system with $\theta_0 = 0.2$. The two excitation profiles result in regret values of 8.8839 and 8.3221, respectively.}}
    \label{fig:profile}
    \vspace{-0.5cm}
    
\end{figure}

For the two remaining systems with $\theta_0 = -0.4$ and $-0.6$, decaying Gaussian and immediate Gaussian exploration gave the same regret with both designs $(a)$ and $(b)$ and the exponent $p$ which was picked was $-20$ in both cases. Hence, the optimal decaying Gaussian exploration is very close to an immediate Gaussian exploration. 

\vspace{-0.08cm}\subsection{Analysis of a particular system}
\vspace{-0.081cm}

In this paragraph, we discuss the results for one particular system with $\theta_0 = -0.4$. In Table~\ref{tab:regret_both_methods}, we give the regret $\bar{R}$  obtained with the four explorations tuned following designs $(a)$ and $(b)$. First, we observe that, for each exploration, both designs give regret which are close to each other, showing that $R_{ub}$ is not far from the actual regret $\bar{R}$ and so minimizing $R_{ub}$ is a good practice in order to minimize the regret $\bar{R}$. It is also clear that immediate binary exploration provides a much lower regret than using immediate Gaussian exploration. \textcolor{black}{The reason behind the difference will be studied in future work.}
\textcolor{black}{Moreover, the optimal grid-search exponent $p$ for the decaying Gaussian exploration was chosen as $-20$, diverging from the commonly used value of $-0.5$ in the LQR literature~\cite{Wan:21,jedra2021minimal}. If $p$ is set to the non-optimal value of $-0.5$, the coefficient $c$ determined by design (b) was 0, leading to the lazy excitation. }

\begin{figure}
    \centering
    \includegraphics[width =0.95\linewidth]{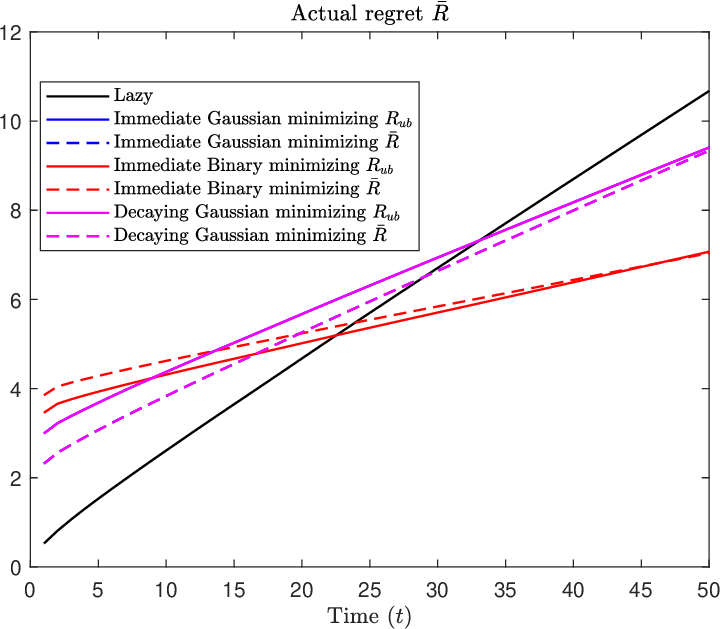}
    \vspace{-0.1cm}
    \caption{Time evolution of the regret with lazy (black solid line), decaying Gaussian (magenta lines), immediate Gaussian (blue lines) and immediate binary (red lines) explorations, tuned with  both designs $(a)$ (solid lines) and $(b)$ (dashed lines). The magenta and blue lines are on top of each other.}
    \label{fig:actual_regret_one_run}
    \vspace{-0.5cm}
    
\end{figure}

In Figure~\ref{fig:actual_regret_one_run}, we depict the time evolution of the  average of the costs $\sum_{k=1}^{t} (\Phi(u^*_k+\alpha_k, \theta_0)-\Phi(u^*_0, \theta_0) )$ obtained from the $1000$ Monte Carlo simulations when the immediate and decaying exploration strategies are optimized for horizon $T = 50$ following designs $(a)$ and $(b)$. 
We notice that lazy exploration is outperformed by immediate exploration already after only half of the design horizon $T$ has elapsed. This illustrates Remark~\ref{rem:interpretation_th1}, i.e. it is advantageous to momentarily degrade the regret with a large exploration at the beginning as it will eventually pay off due to the lower regrets obtained after the exploration (notice that the slopes of the regrets for the immediate explorations are lower than for the lazy exploration). By comparing the time evolution of the regret obtained with design $(a)$ (solid lines) and $(b)$ (dashed lines), we observe that the difference vanishes at the end of the horizon, which justifies our methodology on optimizing $R_{ub}$.

\begin{table}
    \centering
    \caption{Regret $\bar{R}$ obtained with the four explorations tuned following designs $(a)$ and $(b)$.}
    \begin{tabular}{c|c|c}
        $\bar{R}$ & Design $(a)$ &  Design $(b)$ \\ \hline 
        Lazy & 10.676 & 10.676 \\
        Immediate Gaussian & 9.408 & 9.338 \\
        Immediate Binary & 7.070 & 7.039\\
        Decaying Gaussian &  9.408 & 9.338\\
    \end{tabular}
    \label{tab:regret_both_methods}
    \vspace{-0.1cm}
\end{table}

\vspace{-0.3cm}\section{Conclusion}
\label{sec:conclu}

In this work, we analyzed the problem of designing \textcolor{black}{effective exploration} strategies  based on regret minimization in the framework of unconstrained scalar optimization of non-linear static systems. \textcolor{black}{The primary focus of this work is conceptual, aiming to bridge the understanding from linear to non-linear systems.} 
We proposed several approximations to solve this challenging problem and showed that the optimal exploration strategy to minimize regret is either a lazy or an immediate exploration. 
\textcolor{black}{This finding highlights the critical importance of exploration at the onset of the time horizon, significantly simplifying both the design and implementation of exploration strategies.}
This result was supported by a numerical example where we illustrated that minimizing the approximate regret upper bound provides satisfactory performance to minimize the actual regret and that an immediate exploration generated from a zero-mean binary white noise is better than using zero-mean Gaussian white noise. 

In future work, we will conduct a more detailed analysis of the choice of distribution for the exploration signal, such as binary, Gaussian and deterministic. Additionally, we will explore the asymptotic behavior of our exploration strategies.
Finally, further developments are required to obtain a practically useful method to handle the parameter dependency of the approximate regret used for the exploration design.


\bibliographystyle{IEEEtran}
\bibliography{references}

\appendix
\vspace{-0.1cm}

\vspace{-0.08cm}\subsection{Fisher information}

\label{sec:app1}
We will compute the Fisher information for the system in Fig.~\ref{fig:diagram}. 
For simplicity, we assume that the PDF $p_e$ is known ($\theta$ independent) and time invariant. 
The first step is to establish the likelihood function representing the probability of observing $\{u^t, y^t\}$ generated from $y_t = h(u_t,\theta) + e_t$ in iterative optimization with $u_t = u_t^* + \alpha_t$. 
Here the superscript $t$ denotes denotes all past data up to and including $t$. Using several successive chain rules,
\begin{align}
&p(u^t, y^t; \theta) = p(y_t|u^t,y^{t-1}; \theta)p(u_t,y^{t-1}; \theta) \nonumber \\
=&p_{\epsilon}\big(y_t-h(u_t,\theta)\big) p(u_t|u^{t-1},y^{t-1};\theta)p(u^{t-1},y^{t-1};\theta)		  \nonumber \\
=&p_{\epsilon}\big(y_t-h(u_t,\theta)\big) p_{\alpha}(u_t-u^*_t)p(u^{t-1},y^{t-1};\theta)
\nonumber \\
=&\cdots = \prod_{s=1}^{t}\big[p_{\epsilon}\big(y_s-h(u_s,\theta)\big)p_{\alpha}(u_s-u^*_s)\big]
\end{align}
\vspace{-0.3cm}

\noindent where the notation $p(a|b)$ refers to the conditional distribution of $a$ after observing $b$, $p_{\epsilon}$ and $p_{\alpha}$ are the probability density functions of the model  residual and the exploration input, respectively. Note that the second term $p_{\alpha}(u_s-u^*_s)$ is independent of $\theta$. The score function is defined as the gradient w.r.t. $\theta$ of the log-likelihood function:
\vspace{-0.6cm}

\begin{align}
\frac{\partial}{\partial \theta}\log p(u^t, y^t; \theta)=&\sum_{s=1}^{t}\frac{\partial}{\partial \theta} \log p_{\epsilon}\big(y_s-h(u_s,\theta)\big) \\
=&-\sum_{s=1}^{t} \frac{p'_{\epsilon}\big(y_s-h(u_s,\theta)\big)}{p_{\epsilon}\big(y_s-h(u_s,\theta)\big)}\frac{\partial h}{\partial \theta} \nonumber
\end{align}where $p'_{\epsilon}(y_s-h(u_s,\theta))$ is the derivative of $p_{\epsilon}(y_s-h(u_s,\theta))$ w.r.t.   $\theta$ using  chain rule, with $p'_{\epsilon}$ being the derivative of $p_{\epsilon}$ w.r.t.  its argument $y_s-h(u_s,\theta)$.
The Fisher information, denoted by $\mathbb{I}_t$, is a measure of the information that the observed data $\{u^t, y^t\}$ provides about the parameter. It is calculated by squaring the score function and taking its expected value w.r.t. the random variables $\{e_t\}$ and\footnote{From a statistical perspective, $\alpha_t$ is an ancillary statistic meaning that it contains no information about $\theta$. Still such a statistic may heavily influence the properties of an estimator and the Fisher information should be conditioned on such a statistic (see, e.g., \cite{Efron:78a,Fraser:04}).} $\{\alpha_t\}$ at the true parameter $\theta_0$. 
This implies that $\epsilon_s=y_s-h(u_s,\theta_0) = e_s$ and so $\mathbb{I}_t$ is given by~\eqref{eq:It_1}-\eqref{eq:It_2}.
\begin{figure*}
\begin{align}
    &\mathbb{I}_t=\mathbb{E}\bigg[\big(\frac{\partial}{\partial \theta}\log p(u^t, y^t; \theta_0)\big)
\big(\frac{\partial}{\partial \theta}\log p(u^t, y^t; \theta_0) \big)^T
\bigg]\label{eq:It_1}\\
&=\mathbb{E}\bigg[ \sum_{s=1}^{t} \sum_{r=1}^{t} 
p_e'(e_s) p_e'(e_r)
\frac{d \log p_e(e_s)}{d e} 
\frac{d \log p_e(e_r)}{d e} 
\frac{ p_e'(e_s)}{ p_e(e_s)} 
\frac{ p_e'(e_r)}{ p_e(e_r)} 
\frac{\partial h}{\partial \theta}\Bigl|_{\substack{\theta=\theta_0\\u_s=u^*_s+\alpha_s}} 
\frac{\partial h}{\partial \theta}\Bigl|_{\substack{\theta=\theta_0\\u_r=u^*_r+\alpha_r}}\bigg] \label{eq:It_2}
\end{align}
\end{figure*} Since $\{e_t\}$ is a zero-mean white noise, we have independence between $e_s$ and $e_r$ for every $r\ne s$. By recalling that $e_t$ is zero-mean and Gaussian with variance $\sigma^2$, we get
\vspace{-0.4cm}

\begin{align}
&\mathbb{I}_t=  \sum_{s=1}^{t} \frac{1}
{\sigma^{2}} \mathbb{E}\bigg[\frac{\partial h}{\partial \theta}\Bigl|^2_{\substack{\theta=\theta_0\\u_s=u^*_s+\alpha_s}}\bigg] \nonumber
&=\mathbb{E}
\bigg[
\sum_{s=1}^{t}	
\frac{1}
{\sigma^{2}}
\frac{\partial h}{\partial \theta}
\Bigl|_{\substack{\theta=\theta_0\\u_s=u^*_s+\alpha_s}}^2
\bigg] 
\end{align}We obtain the second term in the expression~\eqref{eq:Fish_t}.

It should be noted that the expectation is calculated w.r.t. both $\alpha$ and $e$, which serve as external inputs into the system illustrated in Fig.~\ref{fig:diagram}. 
These inputs are instrumental in engendering the stochastic characteristics of the system illustrated in Fig.~\ref{fig:diagram}. 

\vspace{-0.08cm}\subsection{Proof of Theorem~\ref{th1}}
\label{sec:app2}

We start by proving the following Lemma.
\begin{lemma}
	\label{lemma1}
	Consider the problem \eqref{eq:optimization_problem_theorem} and denote with $x^*=[x^*_1, \dots, x^*_{T-1}]$ its solution.
 Assume that the function $i$ is non-negative and monotonically increasing in the domain $[0, \infty)$. Then $x^*_1\geq x^*_2\geq \cdots \geq x^*_{T-1} \geq 0$.
\end{lemma}

\noindent \textit{Proof:} Consider the vector ${x = [x_1,  \dots   ,x_j, x_{j+1}, \dots ,x_{T-1}]}$ and the vector ${\tilde{x} = [x_1,  \dots,   x_{j+1},  x_{j},  \dots, x_{T-1}]}$ built from $x$ by swapping  $x_j$ and $x_{j+1}$. Assume $x_j\geq x_{j+1}$. Denote with $C : \mathbb{R}^{T-1}_+ \to \mathbb{R}_+$ the objective function of Problem \eqref{eq:optimization_problem_theorem}.
By comparing the cost induced by $x$ and $\tilde{x}$ we get that $C(x)-C(\tilde{x})$ equals to
\begin{align}
	\frac{1}{i_0+\sum_{s=1}^{j-1}i(x_{s})+i(x_{j})} 
-
&\frac{1}{i_0+\sum_{s=1}^{j-1}i(x_{s})+i(x_{j+1})} \leq 0 \nonumber
\end{align}
due to the fact that $x_j \geq x_{j+1}$ and that the information function $i$ is monotonically increasing in the domain $[0, \infty)$. This means that $x_j\geq x_{j+1}$ should be kept to obtain the minimized value. 
The result we want to prove follows from the observation that we can obtain an ordered vector through a finite number of swaps of consecutive elements. Thus 	$x^*_1\geq x^*_2\geq \cdots \geq x^*_{T-1}  \geq 0$.\hfill $\blacksquare$

Now we prove Theorem~\ref{th1} based on Lemma~\ref{lemma1}. 
Firstly, we prove that the optimal solution $x^*$ of Problem~\eqref{eq:optimization_problem_theorem} satisfies $x^*_2=\cdots=x^*_{T-1}=0$, which implies that the optimal solution is either a lazy excitation with $x^*_1=0$ or an immediate excitation with $x^*_1>0$.
We will then discuss a sufficient condition for an immediate excitation to be optimal.

 We start with defining the Lagrangian function for Problem~\eqref{eq:optimization_problem_theorem} as follows
\begin{align}
    L(x,\lambda)=
    \sum_{t=1}^{T-1}
    \bigg[ 
    \frac{1}{
    i_0+\sum_{s=1}^{t}i(x_s)
    }
    +x_t
    +\lambda_t(-x_t)
    \bigg]
\end{align}
where the vector $\lambda=[\lambda_1, \dots, \lambda_T]$ consists of KKT multipliers.
Since the Karush-Kuhn-Tucker (KKT) conditions are necessary for optimality, at the solution $x^*$ it holds
\begin{align}
	1-\sum_{t=k}^{T-1}\frac{i'(x^*_k)}{[i_0+\sum_{s=1}^{t}i(x^*_s)]^2}=\lambda_k^*, \, &k=1,\dots, T-1 \ \ \label{eq:kkt_1}\\
	\lambda_k^*\geq 0,\ x_k^*\geq 0,  \,  &k=1,\dots, T-1 \ \ \label{eq:kkt_2}
 \\
	, \ &k=1,\dots, T-1  \\
\\	\lambda^*_kx_k^*= 0,  \,  &k=1,\dots, T-1\ \ \label{eq:slack}
\end{align}
where $\lambda_k^*$ is the optimal KKT multiplier. 
Due to $x^*_1\geq x^*_k$ for $k=2,\dots,T-1$ (from Lemma~\ref{lemma1}) and the convexity of the information function $i$, which has increasing derivative, we obtain
\begin{align}\label{eq:ineq1}
	i'(x^*_k) \leq i'(x^*_1) \  \Rightarrow \ -i'(x^*_k) \geq -i'(x^*_1) 
\end{align}
Also, for $k=2,\dots,T-1$, it holds that
\begin{equation}\label{eq:ineq2}
	\sum_{t=k}^{T-1}\frac{1}{[i_0+\sum_{s=1}^{t}i(x^*_s)]^2}
	<
	\sum_{t=1}^{T-1}\frac{1}{[i_0+\sum_{s=1}^{t}i(x^*_s)]^2}
\end{equation}
since all the terms in the sum on the right-hand side of~\eqref{eq:ineq2} are positive.
From \eqref{eq:ineq1} and \eqref{eq:ineq2} it follows
\begin{equation}\label{eq:ineq3}
	1-\sum_{t=k}^{T-1}\frac{i'(x^*_k)}{[i_0+\sum_{s=1}^{t}i(x^*_s)]^2}
	>
	1-\sum_{t=1}^{T-1}\frac{i'(x^*_1)}{[i_0+\sum_{s=1}^{t}i(x^*_s)]^2}
	\nonumber
\end{equation}
which together with~\eqref{eq:kkt_1} suggests that for the KKT multipliers it holds $\lambda_k^*>\lambda_1^*\geq 0$, for all $k \ge 2$, which together with the slackness complementary condition in~\eqref{eq:slack}, implies $x_2^*=\cdots=x_{T-1}^*=0$.
Thus we proved that the optimal solution $x^*$ is either a lazy or an immediate excitation.\\
Now assume that the solution $x^*$ is a lazy excitation. Then, according to~\eqref{eq:kkt_1} and~\eqref{eq:kkt_2},
it should hold
\begin{equation}\label{eq:lazy_eq}
	1-\sum_{t=1}^{T-1}\frac{i'(0)}{[i_0+t
		\cdot i(0)]^2}\geq 0
\end{equation}
which proves that the violation of~\eqref{eq:lazy_eq}, i.e.
\begin{equation}\label{eq:immed_condition_again}	
	\sum_{t=1}^{T-1}\frac{i'(0)}{[i_0+t\cdot i(0)]^2}> 1
\end{equation}is sufficient for immediate excitation, since we know that the solution must be either immediate or lazy. \hfill $\blacksquare$

\end{document}